\documentclass[a4paper,12pt]{amsart}
\usepackage{amsfonts}
\usepackage{amssymb}
\usepackage{ifthen}
\usepackage{graphicx}
\nonstopmode \numberwithin{equation}{section}

\usepackage[usenames]{color}
\newtheorem{thm}{Theorem}[section]
\newtheorem{lem}{Lemma}[section]
\newtheorem{cor}[thm]{Corollary}
\newtheorem{prop}[thm]{Proposition}

\newtheorem{step}{Step}[section]

\theoremstyle{definition}
\newtheorem{mlem}{Main lemma}[section]
\newtheorem{assertion}{Assertion}[section]
\newtheorem{cl}{Claim}[section]
\newtheorem{ca}{Case}[section]
\newtheorem{sca}{Subcase}[section]
\newtheorem{scl}{Subclaim}[section]
\newtheorem{conj}[thm]{Conjecture}
\newtheorem{fact}{Fact}[section]
\newtheorem{defn}[thm]{Definition}
\newtheorem{op}[thm]{Open Problem}
\newtheorem{ques}[thm]{Question}
\newtheorem{rem}[thm]{Remark}
\newtheorem{exam}[thm]{Example}

\numberwithin{equation}{section}

\newcounter {own}
\def\theown {\thesection       .\arabic{own}}

\newenvironment{pf}[1][]{%
 \vskip 3mm
 \noindent
 \ifthenelse{\equal{#1}{}}%
  {{\slshape Proof. }}%
  {{\slshape #1.} }%
 }%
{\qed\bigskip}

\newcounter{alphabet}
\newcounter{tmp}
\newenvironment{Thm}[1][]{\refstepcounter{alphabet}%
\bigskip%
\noindent%
{\bf Theorem \Alph{alphabet}}%
\ifthenelse{\equal{#1}{}}{}{ (#1)}%
{\bf .} \itshape}{\vskip 8pt}

\makeatletter
\newcommand{\Ref}[1]{\@ifundefined{r@#1}{}{\setcounter{tmp}{\ref{#1}}\Alph{tmp}}}
\makeatother

\newenvironment{Lem}[1][]{\refstepcounter{alphabet}%
\bigskip%
\noindent%
{\bf Lemma \Alph{alphabet}}%
{\bf .} \itshape}{\vskip 8pt}

\newcounter{alphabet2}

\newcommand{\IC}{{\mathbb C}}
\newcommand{\ID}{{\mathbb D}}

\newcommand{\diam}{{\operatorname{diam}}}

\def\be{\begin{equation}}
\def\ee{\end{equation}}

\newcommand{\ben}{\begin{enumerate}}
\newcommand{\een}{\end{enumerate}}

\newcommand{\blem}{\begin{lem}}
\newcommand{\elem}{\end{lem}}
\newcommand{\bthm}{\begin{thm}}
\newcommand{\ethm}{\end{thm}}
\newcommand{\bcor}{\begin{cor}}
\newcommand{\ecor}{\end{cor}}
\newcommand{\beg}{\begin{exam}}
\newcommand{\eeg}{\end{exam}}
\newcommand{\begs}{\begin{examples}}
\newcommand{\eegs}{\end{examples}}
\newcommand{\bdefe}{\begin{defn}}
\newcommand{\edefe}{\end{defn}}
\newcommand{\bprob}{\begin{prob}}
\newcommand{\eprob}{\end{prob}}
\newcommand{\bques}{\begin{ques}}
\newcommand{\eques}{\end{ques}}
\newcommand{\bei}{\begin{itemize}}
\newcommand{\eei}{\end{itemize}}
\newcommand{\bcon}{\begin{conj}}
\newcommand{\econ}{\end{conj}}
\newcommand{\bop}{\begin{op}}
\newcommand{\eop}{\end{op}}

\newcommand{\bas}{\begin{assertion}}
\newcommand{\eas}{\end{assertion}}

\newcommand{\bfa}{\begin{fact}}
\newcommand{\efa}{\end{fact}}

\newcommand{\bca}{\begin{ca}}
\newcommand{\eca}{\end{ca}}

\newcommand{\bst}{\begin{step}}
\newcommand{\est}{\end{step}}

\newcommand{\bsca}{\begin{sca}}
\newcommand{\esca}{\end{sca}}

\newcommand{\bcl}{\begin{cl}}
\newcommand{\ecl}{\end{cl}}

\newcommand{\bmlem}{\begin{mlem}}
\newcommand{\emlem}{\end{mlem}}

\newcommand{\bscl}{\begin{scl}}
\newcommand{\escl}{\end{scl}}

\newcommand{\bcons}{\begin{conjs}}
\newcommand{\econs}{\end{conjs}}

\newcommand{\bprop}{\begin{prop}}
\newcommand{\eprop}{\end{prop}}

\newcommand{\br}{\begin{rem}}
\newcommand{\er}{\end{rem}}
\newcommand{\brs}{\begin{rems}}
\newcommand{\ers}{\end{rems}}
\newcommand{\bo}{\begin{obser}}
\newcommand{\eo}{\end{obser}}
\newcommand{\bos}{\begin{obsers}}
\newcommand{\eos}{\end{obsers}}
\newcommand{\bpf}{\begin{pf}}
\newcommand{\epf}{\end{pf}}
\newcommand{\ba}{\begin{array}}
\newcommand{\ea}{\end{array}}
\newcommand{\beq}{\begin{eqnarray}}
\newcommand{\beqq}{\begin{eqnarray*}}
\newcommand{\eeq}{\end{eqnarray}}
\newcommand{\eeqq}{\end{eqnarray*}}

\newcounter{minutes}
\divide\time by 60
\newcounter{hours}
\multiply\time by 60 \addtocounter{minutes}{-\time}

\begin{document}

\bibliographystyle{amsplain}
\title []
{The boundary behaviour of $K$-quasiconformal harmonic mappings}

\def\thefootnote{}
\footnotetext{ \texttt{\tiny File:~\jobname .tex,
          printed: \number\day-\number\month-\number\year,
          \thehours.\ifnum\theminutes<10{0}\fi\theminutes}
} \makeatletter\def\thefootnote{\@arabic\c@footnote}\makeatother

\author{Shaolin Chen}
 \address{Sh. Chen, College of Mathematics and
Statistics, Hengyang Normal University, Hengyang, Hunan 421008,
People's Republic of China.} \email{mathechen@126.com}

\author{Saminathan Ponnusamy 
}
\address{S. Ponnusamy, Department of Mathematics,
Indian Institute of Technology Madras, Chennai-600 036, India. }
\email{samy@iitm.ac.in}



\begin{abstract}
In this article,
we first discuss the Lipschitz characteristic and the linear measure distortion of
$K$-quasiconformal harmonic mappings. Then we  give some
characterizations of the radial John disks with the help of
Pre-Schwarzian of  harmonic mappings.
\end{abstract}

\maketitle \pagestyle{myheadings} \markboth{ Sh. Chen and S.
Ponnusamy}{The boundary behaviour of $K$-quasiconformal harmonic
mappings}

\section{Preliminaries and the statement of main results}\label{csw-sec1}


The purpose of this article is to continue our investigations of the boundary behavior of $K$-quasiconformal harmonic mappings,
using the Lipschitz continuity and Pre-Scwarzian derivative defined in \cite{HM}.

\subsection{Notation}
Let $\mathbb{D}$ be the open unit disk in the complex plane $\IC$.
For a sense-preserving harmonic mapping $f=h+\overline{g}$ of $\ID$,  where $h$ and $g$ are analytic in $\ID$,
the Jacobian of $f$ is given by
$J_{f}(z)=  |h'(z)|^2-|g'(z)|^2$ and $\omega (z)=g'(z)/h'(z)$ denotes the dilatation of $f$. Also, we let
$$\|D_{f}\|=|f_{z}|+|f_{\overline{z}}| ~\mbox{ and }~ l(D_{f})=\big|
|f_{z}|-|f_{\overline{z}}|\big |,
$$
where $f_{z}$ and $f_{\overline{z}}$ are the usual partial
derivatives. For $z\in\mathbb{D}$, let
$$B(z)=\{\zeta:\,|z|\leq|\zeta|<1,~|\arg z-\arg \zeta|\leq\pi(1-|z|)\},
$$
and
$$I(z)=\{\zeta\in \partial\mathbb{D}:\,|\arg z-\arg \zeta|\leq\pi(1-|z|)\}.
$$
Let $d_{\Omega}(z)$ be the Euclidean distance from $z$ to the boundary
$\partial \Omega$ of $\Omega$. If $\Omega =\mathbb{D}$, then we set
$d(z):=d_{\mathbb{D}}(z)$. Throughout of this paper, we use the
symbol $C$ to denote the various positive constants, whose value may
vary from one occurrence to another.

\subsection{Preliminaries and Definitions}

\begin{defn}\label{CP-de-1}
A bounded simply connected plane domain $G$ is called a {\it
$c$-John disk} for $c\geq1$ with {\it John center} $w_{0}\in G$ if
for each $w_{1}\in G$ there is a rectifiable arc $\gamma$, called a
{\it John curve}, in $G$ with end points $w_{1}$ and $w_{0}$ such
that
$$\sigma_{\ell}(w)\leq cd_{G}(w)
$$
for all $w$ on $\gamma$, where $\gamma[w_{1},w]$ is the subarc of
$\gamma$ between $w_{1}$ and $w$, and $\sigma_{\ell}(w)$ is the
Euclidean length of $\gamma[w_{1},w]$ (see \cite{COC,KH,John}).
\end{defn}

We can classify $c$-John disk according to some test mappings. More
precisely, if $f$ is a complex-valued and univalent mapping ($f$ is
not necessarily analytic) in $\mathbb{D}$, $G=f(\mathbb{D})$ and,
for $z\in\mathbb{D}$, $\gamma=f([0,z])$ in Definition \ref{CP-de-1},
then we call $c$-John disk  a {\it radial} $c$-John disk, where
$w_{0}=f(0)$ and $w=f(z)$. In particular, if $f$ is univalent and
analytic,  then we call $c$-John disk  a {\it hyperbolic} $c$-John
disk with respect to $f$.  It is well known that any point $w_{0}\in
G$ can be chosen as a John center by modifying the constant $c$ if
necessary. When we do not wish to emphasize the role of $c$, then we
regard the $c$-John disk simply as a John disk in the natural way
(cf. \cite{CP,CP-1,KH,John}).

A sense-preserving homeomorphism $f$ from a domain $\Omega$ onto
$\Omega'$, contained in the {\it Sobolev class}
$W_{loc}^{1,2}(\Omega)$, is said to be a {\it $K$-quasiconformal
mapping} if, for $z\in\Omega$,
$$\|D_{f}(z)\|^{2}\leq K\big | \det D_{f}(z)\big |,~\mbox{i.e.,}~\|D_{f}(z)\|\leq Kl\big(D_{f}(z)\big),
$$
where $K\geq1$ is a constant (cf. \cite{K,LV}).


Let ${\mathcal S}_{H}$ denote the family of sense-preserving planar
harmonic univalent mappings $f=h+\overline{g}$ in $\mathbb{D}$
satisfying the normalization $h(0)=g(0)=h'(0)-1=0$, where $h$ and
$g$ are analytic in $\mathbb{D}$. Recall that $f$ is
sense-preserving in $\mathbb{D}$ if $J_{f}>0 $ in $\mathbb{D}$.
Thus, $f$ is locally univalent and sense-preserving in $\mathbb{D}$
if and only if $J_{f}>0$ in $\mathbb{D}$; or equivalently if $h'\neq
0$ in $\mathbb{D}$ and the dilatation $\omega =g'/h'$ has the
property that $|\omega|<1$ in $\mathbb{D}$ (see
\cite{Clunie-Small-84,Du,Lewy}). The family ${\mathcal S}_{H}$
together with a few other geometric subclasses, originally
investigated in detail by \cite{Clunie-Small-84,Small}, became
instrumental in the study of univalent harmonic mappings (see
\cite{Du}) and has attracted much attention of many function
theorists. If the co-analytic part $g$ is identically zero in the
decomposition of $f=h+\overline{g}$, then the class ${\mathcal
S}_{H}$ reduces to the classical family $\mathcal S$ of all
normalized analytic univalent functions
$h(z)=z+\sum_{n=2}^{\infty}a_{n}z^{n}$ in $\mathbb{D}$. If
${\mathcal S}_H^{0}=\{f=h+\overline{g} \in {\mathcal S}_H: \,g'(0)=0
\} $, then the family ${\mathcal S}_H^{0}$ is both normal and
compact (cf. \cite{Clunie-Small-84, Du}). Denote by ${\mathcal
S}_{H}(K)$ (resp. ${\mathcal S}_{H}^{0}(K)$) if $f\in {\mathcal
S}_{H}$ (resp. ${\mathcal S}_{H}^{0}$) and is a $K$-quasiconformal
harmonic mapping in $\mathbb{D}$, where $K\geq 1$ is a constant.
Also, we denote by ${\mathcal S}_{H}(K,\Omega)$ (resp. ${\mathcal
S}_{H}^{0}(K,\Omega)$) if $f\in {\mathcal S}_{H}(K)$ (resp. $f\in
{\mathcal S}_{H}^{0}(K)$) and $f$ maps $\mathbb{D}$ onto $\Omega$,
where $\Omega$ is a subdomain of $\mathbb{C}$.

\subsection{Statement of Main results}
We now state our first main result  which concerns the Lipschitz
continuity on $K$-quasiconformal harmonic mappings of $\mathbb{D}$
onto a radial John disk.

\begin{thm}\label{thm-1}
Let $f\in {\mathcal S}_{H}^{0}(K,\Omega)$, where $\Omega$ is a radial John disk. Then, for $z\in\mathbb{D}$ and
$\zeta_{1},\,\zeta_{2}\in B(z)$, there are constants $\delta\in(0,1)$
and $C>0$ such that
$$|f(\zeta_{1})-f(\zeta_{2})|\leq Cd_{\Omega}(f(z))\left(\frac{|\zeta_{1}-\zeta_{2}|}{1-|z|}\right)^{\delta}.
$$
\end{thm}

We would like to point out that Theorem \ref{thm-1} was established
in \cite[Theorem 4]{CP-1} but with an additional assumption that
$|z|\geq\frac{1}{2}$, and thus, we see now that the condition
``$|z|\geq\frac{1}{2}$" in \cite[Theorem 4]{CP-1} is redundant.
Moreover, by \cite[Lemma 6]{CP-1} and \cite[Inequality (2.3)]{CP-1},
we obtain
$$\frac{1}{16K}\leq d_{\Omega}(f(0))\leq\frac{2K}{1+K}\|D_{f}(0)\|=\frac{2K}{1+K},
$$
where $\Omega=f(\mathbb{D})$, $\|D_{f}(0)\|=|f_{z}(0)|+|f_{\overline{z}}(0)|$ and
$f\in{\mathcal S}_{H}^{0}(K, \Omega)$. Therefore, by letting $z=0$ in Theorem \ref{thm-1},  we get the following result.

\begin{cor}\label{cor-1}
Let $f\in {\mathcal S}_{H}^{0}(K,\Omega)$, where $\Omega$ is a
radial John disk. Then, for all
$\zeta_{1},\,\zeta_{2}\in\mathbb{D}$, there are constants $C>0$ and
$\delta\in(0,1)$ such that
$$|f(\zeta_{1})-f(\zeta_{2})|\leq C|\zeta_{1}-\zeta_{2}|^{\delta}.
$$
\end{cor}

Our next result establishes the linear measure distortion on $K$-quasiconformal mappings of $\mathbb{D}$ into a
radial John disk.

\begin{thm}\label{thm-2}
Let $f\in{\mathcal S}_H^{0}(K,\Omega)$, where $\Omega$ is a radial
John disk. Then, for all $z_{1},\,z_{2}\in\mathbb{D}$ with
$|z_{2}|\leq|z_{1}|$, there are constant $C>0$ and $\delta\in(0,1)$
such that
$$\frac{\diam \big(f(B(z_{1}))\big)}{\diam \big(f(B(z_{2}))\big)}\leq C\left(\frac{\ell(I(z_{1}))}{\ell(I(z_{2}))}\right)^{\delta}.
$$
\end{thm}

The Pre-Schwarzian derivative $P_{f}$ of a sense-preserving harmonic mapping $f=h+\overline{g}$ in $\mathbb{D}$  is defined by
$$P_f=(\log J_f)_z=\frac{h''\overline{h'}-g''\overline{g'}}{|h'|^2-|g'|^2}= T_h -\frac{\omega'\overline{\omega}}{1-|\omega|^{2}},
$$
where $\omega  =g'/h'$, and $T_h= h''/h'$ denotes the Pre-Schwarzian
of a locally univalent analytic function $h$ in $\mathbb{D}$. See
\cite{CDO,HM,LiPo18} for recent investigations on Pre-Schwarzian
derivatives of harmonic mappings.

Ahlfors and Weill \cite{AW} and, Becker and Pommerenke \cite{BP}
characterized quasidisks by using the Pre-Schwarzian of analytic
functions. On the basis of the works of Chuaqui, et al. \cite{COC},
Kari Hag and Per Hag \cite{KH} discussed relationships between
John disks and the Pre-Schwarzian of analytic functions. By
analogy with \cite[Theorem 4]{COC} and \cite[Theorem 3.7]{KH},
the present authors in \cite[Theorem 5]{CP} showed that if
$f\in {\mathcal S}_{H}^{0}(K,\Omega)$ such that
$$\limsup_{|z|\rightarrow1^{-}}\left\{(1-|z|^{2}){\rm Re}\big(zP_{f}(z)\big)\right\}<1,
$$
then $\Omega$ is a radial John disk. Our final result improves this result in the following form.

\begin{thm}\label{thm-4.0}
Let $f=h+\overline{g}\in {\mathcal S}_{H}^{0}(K,\Omega)$ and $\omega
 =g'/h'$. Then the following statements are true.
\begin{enumerate}
\item[{\rm (a)}] If
\be\label{thm-4}\limsup_{|z|\rightarrow1^{-}}\left\{(1-|z|^{2}){\rm
Re}\big(zP_{f}(z)\big)\right\}<1+k,
\ee
then $\Omega$ is a radial John disk, where $k=\frac{K-1}{K+1}\leq\frac{1}{2}.$

\item[{\rm (b)}] If $h$ is univalent in $\mathbb{D}$ and satisfies
\be\label{chen-1}\limsup_{|z|\rightarrow1^{-}}\left\{(1-|z|^{2})\left|P_{f}(z)
+\frac{\omega'(z)\overline{\omega(z)}}{1-|\omega(z)|^{2}}\right|\right\}<2,
\ee
then $\Omega$ is a radial John disk.
\end{enumerate}
\end{thm}

\begin{cor}\label{cor-1.6}
Let $f=h+\overline{g}\in {\mathcal S}_{H}^{0}(K,\Omega)$ and $\omega
 =g'/h'$. If $h$ is univalent in $\mathbb{D}$ and satisfies
\be\label{chen-cp}\sup_{z\in\mathbb{D}}\left\{(1-|z|^{2})\left|P_{f}(z)
+\frac{\omega'(z)\overline{\omega(z)}}{1-|\omega(z)|^{2}}\right|\right\}<2,
\ee then $\Omega$ is a radial John disk. The constant $2$ is the
best possible.
\end{cor}

The proofs of Theorems \ref{thm-1}, \ref{thm-2}, \ref{thm-4.0} and
Corollary \ref{cor-1.6} will be presented in Section \ref{csw-sec2}.






\section{Proofs of the main results}\label{csw-sec2}
The hyperbolic plane is the unit disk $\ID$ with the hyperbolic
metric
$$\lambda_{\mathbb{D}}(z)|dz|=\frac{|dz|}{1-|z|^{2}}
$$
which is indeed a mapping which associates to each smooth curves
$\gamma$ in $\mathbb{D}$ its hyperbolic length $\ell _{\ID}(\gamma)$
defined by
$$\ell _{\ID}(\gamma)=\int_{\gamma}\lambda_{\mathbb{D}}(z)|dz|=\int_a^b\frac{|z'(t)|}{1-|z(t)|^{2}}\,dt,
$$
where $\gamma$ is parameterized by $z(t)$, $a\leq t\leq b$. The
hyperbolic distance (or {\it Poincar\'e distance})
$\lambda_{\mathbb{D}}(z_{1}, z_{2})$ between points $z_{1}$ and
$z_{2}$ in  $\mathbb{D}$ is then defined by
$$\lambda_{\mathbb{D}}(z_{1}, z_{2})=\inf_{\gamma} \ell _{\ID}(\gamma) =
\tanh^{-1}\left|\frac{z_{1}-z_{2}}{1-\overline{z}_{1}z_{2}}\right|,
$$
where the infimum is taken over all smooth curves $\gamma$ in
$\mathbb{D}$ that joins $z_{1}$ to $z_{2}$ in $\mathbb{D}$ (cf.
\cite{Po1}).

\begin{Lem}{\rm (\cite[Lemma 1]{CP})} \label{Lem-CP-1}
Let $f\in{\mathcal S}_{H}$. Then for $z_{1},~z_{2}\in\mathbb{D}$,
$$\frac{1}{2}\|D_{f}(z_{1})\|e^{-(1+\alpha)\lambda_{\mathbb{D}}(z_{1},z_{2})}\leq\|D_{f}(z_{2})\|\leq
2\|D_{f}(z_{1})\|e^{(1+\alpha)\lambda_{\mathbb{D}}(z_{1},z_{2})},
$$
where  $ \alpha :=\sup_{f\in{\mathcal
S}_{H}}\frac{|h''(0)|}{2}<+\infty$.
\end{Lem}

We remark that $2\leq\alpha<+\infty$, but the sharp value of
$\alpha$ is still unknown (cf. \cite{Clunie-Small-84, Du,Small}).

\begin{Thm}{\rm (\cite[Theorem 2]{CP})}\label{Thm-B}
Let $f\in {\mathcal S}_{H}^{0}(K,\Omega)$, where
$\Omega:=f(\mathbb{D})$ is a bounded domain. Then the following
conditions are equivalent:
\begin{enumerate}
\item[{\rm (1)}]  $\Omega$ is a radial John disk;

\item[{\rm (2)}] There is a positive constant $C$ such that for all $ z\in\mathbb{D}$,
$$\diam f(B(z))\leq C d_{\Omega}(f(z));
$$
\item[{\rm (3)}] There are  constants $C>0$  and $\delta\in(0,1)$ such
that for all $ z\in\mathbb{D}$ and $\zeta\in B(z)$,
$$\|D_{f}(\zeta)\|\leq C\|D_{f}(z)\|\left(\frac{1-|\zeta|}{1-|z|}\right)^{\delta-1}.
$$
\end{enumerate}
\end{Thm}

\begin{Lem}{\rm (\cite[Lemma 2]{CP})} \label{Lem-CP}
Let $a_{1}, a_{2}$ and $a_{3}$ be positive constants and let
$0<|z_{0}|=1-\delta_{0}$, where $\delta_{0}\in(0,1)$. If
$f=h+\overline{g}\in{\mathcal S}_{H}$,
$0\leq1-a_{2}\delta_{0}\leq|z|\leq1-a_{1}\delta_{0}$ and $|\arg z-
\arg z_{0}|\leq a_{3}\delta_{0}$, then
$$\frac{1}{M(a_{1},a_{2},a_{3})}\|D_{f}(z_{0})\|\leq\|D_{f}(z)\|\leq
M(a_{1},a_{2},a_{3})\|D_{f}(z_{0})\|,
$$
where
$M(a_{1},a_{2},a_{3})=2e^{(1+\alpha)\left(a_{3}+\frac{1}{2}\log\frac{2a_{2}-a_{1}}{a_{1}}\right)}$
and $ \alpha :=\sup_{f\in{\mathcal S}_{H}}\frac{|h''(0)|}{2}. $
\end{Lem}

\begin{Lem}{\rm (\cite[Lemma 6]{CP-1})}\label{Lem-ch-1}
If $f\in{\mathcal S}_{H}(K)$ and $\Omega=f(\mathbb{D})$,  then for $z\in\mathbb{D}$,
$$d_{\Omega}(f(z))\geq\frac{\|D_{f}(z)\|(1-|z|^{2})}{16K}.
$$
\end{Lem}

\subsection{Proof of Theorem \ref{thm-1}} Let $z=re^{i\theta}$,
$\mu=|\zeta_{1}-\zeta_{2}|$ and
$\zeta_{j}=r_{j}e^{i\theta_{j}}~(j=1,2)$ with $r_{1}\leq r_{2}$.

{\rm $\mathbf{Case~ I.}$} If $\rho=1-2\mu<r$, then
$\frac{2\mu}{1-r}>1$ and, by Theorem \Ref{Thm-B}(2), we see that
there is a positive constant $C$ such that
\be \label{c-1}
|f(\zeta_{1})-f(\zeta_{2})|\leq \diam(B(z))\leq Cd_{\Omega}(f(z))
\leq
2^{\delta}Cd_{\Omega}(f(z))\left(\frac{|\zeta_{1}-\zeta_{2}|}{1-|z|}\right)^{\delta}.
\ee

{\rm $\mathbf{Case~ II.}$} Suppose that $\rho=1-2\mu\geq r$ and
$r_{1}<\rho$. In this case, for $|\zeta-\zeta_{1}|\leq\mu$, we have
$$\frac{|\zeta-\zeta_{1}|}{|1-\overline{\zeta}\zeta_{1}|}\leq\frac{\mu}{1-r_{1}}<\frac{\mu}{1-\rho}=\frac{1}{2},
$$
which implies that
\be \label{c-2}
\lambda_{\mathbb{D}}(\zeta, \zeta_{1})=
\tanh^{-1}\left|\frac{\zeta-\zeta_{1}}{1-\overline{\zeta}\zeta_{1}}\right|\leq \tanh^{-1}\left( \frac{1}{2}\right )=\frac{1}{2}\log3,
\ee
where $\lambda_{\mathbb{D}}(z_{1}, z_{2}) $ denotes the hyperbolic distance (or {\it Poincar\'e distance})  between points
$z_{1}$ and $z_{2}$ in  $\mathbb{D}$ given by
$$\lambda_{\mathbb{D}}(z_{1}, z_{2}) =
\tanh^{-1}\left|\frac{z_{1}-z_{2}}{1-\overline{z}_{1}z_{2}}\right|.
$$

It follows from (\ref{c-2}) and Lemma \Ref{Lem-CP-1} that there is a
positive constant $C$ such that
\be \label{c-3} \|D_{f}(\zeta)\|\leq
C\|D_{f}(\zeta_{1})\|,
\ee
where $|\zeta-\zeta_{1}|\leq\mu$. By (\ref{c-3}), it follows that
\beq\label{c-4}
|f(\zeta_{1})-f(\zeta_{2})|&\leq&\int_{[\zeta_{1},\zeta_{2}]}\|D_{f}(\zeta)\||d\zeta|\leq
C|\zeta_{1}-\zeta_{2}|\|D_{f}(\zeta_{1})\|,
\eeq
where
$[\zeta_{1},\zeta_{2}]$ is  the line segment from $\zeta_{1}$ to $\zeta_{2}$.

By Theorem \Ref{Thm-B}(3)  and Lemma \Ref{Lem-ch-1}, there are
constants $C>0$ and $\delta\in(0,1)$ such that
\beqq
\|D_{f}(\zeta_{1})\|&\leq&
C\|D_{f}(z)\|\left(\frac{1-r_{1}}{1-r}\right)^{\delta-1}\leq
C\|D_{f}(z)\|\left(\frac{1-\rho}{1-r}\right)^{\delta-1}\\ &\leq&
16KCd_{\Omega}(f(z))\frac{(1-\rho)^{\delta-1}}{(1-r)^{\delta}},
\eeqq
which, together with (\ref{c-4}), implies that there is a
positive constant $C$ such that
\be\label{c-5}
|f(\zeta_{1})-f(\zeta_{2})|\leq
Cd_{\Omega}(f(z))\frac{(1-\rho)^{\delta-1}}{(1-r)^{\delta}}|\zeta_{1}-\zeta_{2}|\leq
2^{\delta-1}C\left(\frac{|\zeta_{1}-\zeta_{2}|}{1-|z|}\right)^{\delta}.
\ee

{\rm $\mathbf{Case~ III.}$} Suppose that $r\leq\rho=1-2\mu\leq
r_{1}.$ Then, by  Theorem \Ref{Thm-B}(3)  and Lemma \Ref{Lem-ch-1},
we see that there are constants $C>0$ and $\delta\in(0,1)$ such that
\beq\label{c-6}
|f(\zeta_{1})-f(\rho e^{i\theta_{1}})|
 &\leq&\int_{\rho}^{r_{1}}\|D_{f}(te^{i\theta_{1}})\|dt\\
\nonumber &\leq &
C\|D_{f}(z)\|\int_{\rho}^{r_{1}}\left(\frac{1-t}{1-r}\right)^{\delta-1}dt\\
\nonumber
  &=&\frac{C}{\delta}\frac{\|D_{f}(z)\|}{(1-r)^{\delta-1}}\left[(1-\rho)^{\delta}-(1-r_{1})^{\delta}\right]\\
 \nonumber&\leq&\frac{C}{\delta}\frac{\|D_{f}(z)\|(1-\rho)^{\delta}}{(1-r)^{\delta-1}}\\ \nonumber
 &\leq &
 \frac{16KC}{2^{\delta}\delta}d_{\Omega}(f(z))\left(\frac{|\zeta_{1}-\zeta_{2}|}{1-r}\right)^{\delta}~\mbox{ (by Lemma \Ref{Lem-ch-1})}
\eeq
and
\beq\label{c-7}
 |f(\zeta_{2})-f(\rho
 e^{i\theta_{2}})|&\leq&\int_{\rho}^{r_{2}}\|D_{f}(te^{i\theta_{2}})\|dt\\ \nonumber
 &\leq&
 C\|D_{f}(z)\|\int_{\rho}^{r_{2}}\left(\frac{1-t}{1-r}\right)^{\delta-1}dt\\
 \nonumber
 &=&\frac{C}{\delta}\frac{\|D_{f}(z)\|}{(1-r)^{\delta-1}}\left[(1-\rho)^{\delta}-(1-r_{1})^{\delta}\right]\\
 \nonumber&\leq&\frac{C}{\delta}\frac{\|D_{f}(z)\|(1-\rho)^{\delta}}{(1-r)^{\delta-1}}\\
 \nonumber&\leq&
 \frac{16KC}{2^{\delta}\delta}d_{\Omega}(f(z))\left(\frac{|\zeta_{1}-\zeta_{2}|}{1-r}\right)^{\delta}.
\eeq


Let $\gamma$ be the smaller subarc of $\partial\mathbb{D}_{\rho}$
between $\rho e^{i\theta_{1}}$ and $\rho e^{i\theta_{2}}$. Then
$|\theta_{1}-\theta_{2}|\leq\pi$ and, since
\begin{eqnarray*}
|\zeta_{1}-\zeta_{2}| &=&
\sqrt{(r_{1}-r_{2})^{2}+4r_{1}r_{2}\sin^{2}\left(\frac{\theta_{1}-\theta_{2}}{2}\right)}\\
&\geq&
2\sqrt{r_{1}r_{2}}\left|\sin\frac{\theta_{1}-\theta_{2}}{2}\right|\\
&\geq&\frac{2\rho|\theta_{1}-\theta_{2}|}{\pi},
\end{eqnarray*}
we see that
\be\label{c-9}
\ell(\gamma) =\rho|\theta_{1}-\theta_{2}|\leq \frac{\pi}{2}\mu.
\ee
Hence,   we get
\beq\label{c-10}
|f(\rho  e^{i\theta_{1}})-f(\rho e^{i\theta_{2}})|
&\leq&\int_{\gamma}\rho\|D_{f}(\rho e^{i\tau})\|d\tau\\
\nonumber &\leq & C\int_{\gamma}\|D_{f}(z)\|\left(\frac{1-\rho}{1-r}\right)^{\delta-1}d\tau ~\mbox{ (by Theorem \Ref{Thm-B}(3))}\\
\nonumber&=&C\|D_{f}(z)\|\left(\frac{1-\rho}{1-r}\right)^{\delta-1}\ell(\gamma)\\
\nonumber &\leq&\frac{\pi}{4}C\|D_{f}(z)\|\frac{(1-\rho)^{\delta}}{(1-r)^{\delta-1}} ~\mbox{ (by (\ref{c-9}))}\\
\nonumber & \leq &4\pi KCd_{\Omega}(f(z))\left(\frac{1-\rho}{1-r}\right)^{\delta} ~\mbox{ (by Lemma \Ref{Lem-ch-1})}\\
 \nonumber&=&2^{2+\delta}\pi KCd_{\Omega}(f(z))\left(\frac{|\zeta_{1}-\zeta_{2}|}{1-r}\right)^{\delta}.
 \eeq

%

Therefore, by (\ref{c-6}), (\ref{c-7}) and  (\ref{c-10}), we
conclude that there is a positive constant $C$ such that
\beqq
|f(\zeta_{1})-f(\zeta_{2})|&\leq&|f(\zeta_{1})-f(\rho
 e^{i\theta_{1}})|+|f(\zeta_{2})-f(\rho
 e^{i\theta_{2}})|+ |f(\rho  e^{i\theta_{1}})-f(\rho
 e^{i\theta_{2}})|\\
 &\leq&Cd_{\Omega}(f(z))\left(\frac{|\zeta_{1}-\zeta_{2}|}{1-r}\right)^{\delta}.
 \eeqq
The proof of the theorem is complete. \hfill $\Box$

\begin{Thm}{\rm (\cite[Theorem 1]{CP})}\label{Thm-CP-1}
Let $f\in {\mathcal S}_{H}^{0}(K,\Omega)$, where
$\Omega:=f(\mathbb{D})$ is a bounded domain. Then $\Omega$ is a
radial John disk if and only if there are constants $M(K)>0$ and
$\delta\in(0,1)$ such that for each $\zeta\in\partial\mathbb{D}$ and
for $0\leq r\leq\rho<1$,
$$\|D_{f}(\rho\zeta)\|\leq M(K)\|D_{f}(r\zeta)\|\left(\frac{1-\rho}{1-r}\right)^{\delta-1}.
$$
\end{Thm}

\subsection{Proof of Theorem \ref{thm-2}} 
Let $f\in{\mathcal S}_H^{0}(K,\Omega)$, where $\Omega$ is a radial
John disk. Suppose that $z_{1}=re^{i\theta}$ and
$r_{1}e^{i\theta_{1}},r_{2}e^{i\theta_{2}}\in B(z_{1})$ with
$r_{2}\leq r_{1}$. Then, by Theorem \Ref{Thm-B}(3), we see that
there are  positive constants $C$ and $\delta\in(0,1)$ such that
\beq\label{c-eq11}
|f(r_{1}e^{i\theta_{1}})-f(re^{i\theta_{1}})|&\leq&\int_{r}^{r_{1}}\|D_{f}(\rho
e^{i\theta_{1}})\|d\rho\leq
C\int_{r}^{r_{1}}\|D_{f}(z_{1})\|\left(\frac{1-\rho}{1-r}\right)^{\delta-1}d\rho\\
\nonumber
&=&\frac{C}{\delta}\frac{\|D_{f}(z_{1})\|}{(1-r)^{\delta-1}}\left[(1-r)^{\delta}-(1-r_{1})^{\delta}\right]\\
\nonumber&\leq&\frac{C}{\delta} \|D_{f}(z_{1})\|(1-r).
\eeq
Similarly, we have
\beq\label{c-eq12}
|f(r_{2}e^{i\theta_{2}})-f(re^{i\theta_{2}})|
 \leq\frac{C}{\delta} \|D_{f}(z_{1})\|(1-r).
\eeq
Let $\gamma$ be the smaller subarc of $\partial\mathbb{D}_{r}$
between $re^{i\theta_{1}}$ and $re^{i\theta_{2}}$. Since
$r_{1}e^{i\theta_{1}},r_{2}e^{i\theta_{2}}\in B(z_{1})$, we see that
\be\label{c-eq13}
|\theta_{1}-\theta_{2}|\leq|\theta_{1}-\theta|+|\theta-\theta_{1}|\leq2\pi(1-r).
\ee
It follows from (\ref{c-eq13}) and Theorem \Ref{Thm-B}(3) that
\beq\label{c-eq14}
|f(re^{i\theta_{1}})-f(re^{i\theta_{2}})|&\leq&r\int_{\gamma}\|D_{f}(r
e^{i\eta})\|d\eta\leq C\int_{\gamma}\|D_{f}(r e^{i\theta})\|d\eta\\
\nonumber &=& Cr\|D_{f}(r e^{i\theta})\||\theta_{1}-\theta_{2}| \leq
2C\pi(1-r)\|D_{f}(r e^{i\theta})\|.
\eeq
Combining (\ref{c-eq11}), (\ref{c-eq12}) and (\ref{c-eq14}) shows that
\beqq
|f(r_{1}e^{i\theta_{1}})-f(r_{2}e^{i\theta_{2}})|&\leq&|f(r_{1}e^{i\theta_{1}})-f(re^{i\theta_{1}})|+
|f(r_{2}e^{i\theta_{2}})-f(re^{i\theta_{2}})| \\&&
+|f(re^{i\theta_{1}})-f(re^{i\theta_{2}})|\\
&\leq&\left(2\pi C+\frac{2C}{\delta}\right)(1-r)\|D_{f}(r
e^{i\theta})\|,
\eeqq
which implies that there is a positive constant $C$ such that
\be\label{c-eq15}
\diam(B(z_{1}))\leq C(1-|z_{1}|^{2})\|D_{f}(z_{1})\|.
\ee

It follows from  Theorem \Ref{Thm-B}(3), Lemma \Ref{Lem-ch-1} and
\cite[Inequality (2.3)]{CP} that there is a positive constant $C$
such that
\be\label{c-eq16}\diam(B(z_{2}))\geq
d_{\Omega}(f(z_{2}))\geq\frac{\|D_{f}(z_{2})\|(1-|z_{2}|^{2})}{16K}.
\ee

By (\ref{c-eq15}), (\ref{c-eq16}) and Theorem \Ref{Thm-CP-1}, we
conclude that there are constants $M(K)>0$ and $\delta\in(0,1)$ such
that
\beqq \frac{\diam f(B(z_{1}))}{\diam f(B(z_{2}))}&\leq&
16KC\frac{\|D_{f}(z_{1})\|(1-|z_{1}|^{2})}{\|D_{f}(z_{2})\|(1-|z_{2}|^{2})}\\
&\leq& 32M(K)KC\left(\frac{1-|z_{1}|}{1-|z_{2}|}\right)^{\delta}\\
&=&32M(K)KC\left(\frac{\ell(I(z_{1}))}{\ell(I(z_{2}))}\right)^{\delta},
\eeqq
 which completes the proof. \hfill $\Box$

\subsection{Proof of Theorem \ref{thm-4.0}} We first prove (a).
 It follows from (\ref{thm-4}) that there is a $\nu\in(0,1+k)$ and
$r_{0}\in(0,1)$ such that, for $r_{0}\leq\eta<1,$
\be\label{eq-12a}
\frac{\nu}{1-\eta^{2}}\geq{\rm Re}\big(\zeta P_{f}(\eta\zeta)\big)=\mbox{Re}\left(\frac{\zeta h''(\eta\zeta)}{h'(\eta\zeta)}\right)-
\mbox{Re}\left(\frac{\zeta\omega'(\eta\zeta)\overline{\omega(\eta\zeta)}}{1-|\omega(\eta\zeta)|^{2}}\right),
\ee
where $\zeta\in\partial\mathbb{D}.$ By Schwarz-Pick's lemma, we obtain
\be\label{eq-12}
|\omega'(\eta\zeta)|\leq\frac{1-|\omega(\eta\zeta)|^{2}}{1-\eta^{2}}
\ee
and, since $f$ is a $K$-quasiconformal harmonic mapping, we see that,
\be\label{eqr12.0}
|\omega(z)|\leq k=\frac{K-1}{K+1}, \quad z\in\mathbb{D}.
\ee
Thus,  by (\ref{eq-12}) and (\ref{eqr12.0}),  \eqref{eq-12a} gives
\beqq
\nonumber\mbox{Re}\left(\frac{\zeta
h''(\eta\zeta)}{h'(\eta\zeta)}\right)&\leq&
\mbox{Re}\left(\frac{\zeta\omega'(\eta\zeta)\overline{\omega(\eta\zeta)}}{1-|\omega(\eta\zeta)|^{2}}\right)+\frac{\nu}{1-\eta^{2}}
\\
\nonumber&\leq&\frac{|\omega'(\eta\zeta)|\, |\overline{\omega(\eta\zeta)|}}{1-|\omega(\eta\zeta)|^{2}}+\frac{\nu}{1-\eta^{2}}\\
\nonumber&\leq&\frac{\nu+k}{1-\eta^{2}}.
\eeqq
Choosing $\lambda\in(0,k+1-\nu)$, there is an $r_{1}\in[r_{0},1)$ such that
\be\label{eq-y1} \mbox{Re}\left(\frac{\zeta
h''(\eta\zeta)}{h'(\eta\zeta)}\right)<\frac{2\eta-(\lambda+1-2k)}{1-\eta^{2}}
~\mbox{ for all $\zeta\in\partial\mathbb{D}$},
\ee
when $\eta\in[r_{1},1)$. For $0\leq r_{1}\leq r\leq\rho<1$, by (\ref{eq-y1}), we find that
\beq \nonumber
\log\left[\frac{(1-\rho^{2})|h'(\rho\zeta)|}{(1-r^{2})|h'(r\zeta)|}\right]&=&\int_{r}^{\rho}\bigg[\mbox{Re}\left(\frac{\zeta
h''(\eta\zeta)}{h'(\eta\zeta)}\right)-\frac{2\eta}{1-\eta^{2}}\bigg]d\eta\\
\nonumber
&<&-2\left(\frac{\lambda+1}{2}-k\right)\int_{r}^{\rho}\frac{d\eta}{1-\eta^{2}}\\
\nonumber&=&-\left(\frac{\lambda+1}{2}-k\right)\log\left(\frac{1+\rho}{1+r}\cdot\frac{1-r}{1-\rho}\right),
\eeq
which implies that
\be\label{eq-y2}
\left|\frac{h'(\rho\zeta)}{h'(r\zeta)}\right|<\left(\frac{1+r}{1+\rho}\right)^{\frac{1+\lambda}{2}-k+1}
\left(\frac{1-\rho}{1-r}\right)^{\frac{\lambda+1}{2}-1-k}\leq\left(\frac{1-\rho}{1-r}\right)^{\frac{\lambda+1}{2}-1-k}.
\ee
By (\ref{eq-y2}), we get
\beq
\nonumber \|D_{f}(\rho\zeta)\|&\leq&\frac{2K}{1+K}|h'(\rho\zeta)|<\frac{2K}{1+K}
|h'(r\zeta)|\left(\frac{1-\rho}{1-r}\right)^{\frac{\lambda+1}{2}-1-k}\\
\nonumber &\leq&\frac{2K}{1+K}
\|D_{f}(r\zeta)\|\left(\frac{1-\rho}{1-r}\right)^{\frac{\lambda+1}{2}-k-1}.
\eeq

Next, we can use the similar approach  as in the proof of
\cite[Theorem 5]{CP} to remove the restriction $r\geq r_{1}$ above.
Hence, for $0\leq r\leq\rho<1$, there is a positive constant $C$
such that
$$\|D_{f}(\rho\zeta)\|\leq
C\|D_{f}(r\zeta)\|\left(\frac{1-\rho}{1-r}\right)^{\left(\frac{\lambda+1}{2}-k\right)-1},
$$
which, together with Theorem \Ref{Thm-CP-1}, implies $\Omega$ is a
radial John disk.

Now we prove the part of (b). Let $f=h+\overline{g}\in {\mathcal
S}_{H}^{0}(K,\Omega)$ satisfy (\ref{chen-1}), where $h$ is
univalent in $\mathbb{D}$. Then
\beqq
\limsup_{|z|\rightarrow1^{-}}\left\{(1-|z|^{2})\left|\frac{h''(z)}{h'(z)}\right|\right\}=
\limsup_{|z|\rightarrow1^{-}}\left\{(1-|z|^{2})\left|P_{f}(z)+\frac{\omega'(z)\overline{\omega(z)}}{1-|\omega(z)|^{2}}\right|\right\}<2,
\eeqq
which implies that
\be\label{eq-chen1}
\limsup_{|z|\rightarrow1^{-}}\left\{(1-|z|^{2})\mbox{Re}\left(z\frac{h''(z)}{h'(z)}\right)\right\}<2.
\ee
It follows from (\ref{eq-chen1}), \cite[Theorem 3.7]{KH} and \cite[Theorem 2.3]{KH} that there are constants $C>0$ and
$\delta\in(0,1)$ such that for each $\zeta\in\partial\mathbb{D}$ and
for $0\leq r\leq\rho<1,$
\be\label{eq-29} |h'(\rho\zeta)|\leq
C|h'(r\zeta)|\left(\frac{1-\rho}{1-r}\right)^{\delta}.
\ee
Since $f$ is a $K$-quasiconformal mapping, we see that
\be\label{eq-30}
\frac{2}{1+K}|h'(z)|\leq\|D_{f}(z)\|\leq\frac{2K}{1+K}|h'(z)|.
\ee
By (\ref{eq-29}) and (\ref{eq-30}), there are constants $C>0$ and
$\delta\in(0,1)$ such that for each $\zeta\in\partial\mathbb{D}$ and
for $0\leq r\leq\rho<1,$
\beqq
\frac{K+1}{2K}\|D_{f}(\rho\zeta)\|\leq|h'(\rho\zeta)|\leq
C|h'(r\zeta)|\left(\frac{1-\rho}{1-r}\right)^{\delta}\leq\frac{(K+1)C}{2}\|D_{f}(r\zeta)\|\left(\frac{1-\rho}{1-r}\right)^{\delta},
\eeqq
which, together with Theorem \Ref{Thm-CP-1}, yields that
$\Omega$ is a radial John disk. The proof of the theorem is
complete. \hfill $\Box$

\subsection{Proof of Corollary \ref{cor-1.6}}
By the assumption, we have

\beqq
\sup_{z\in\mathbb{D}}\left\{(1-|z|^{2})\left|\frac{h''(z)}{h'(z)}\right|\right\}=
\sup_{z\in\mathbb{D}}\left\{(1-|z|^{2})\left|P_{f}(z)+\frac{\omega'(z)\overline{\omega(z)}}{1-|\omega(z)|^{2}}\right|\right\}<2,
\eeqq which implies that \be\label{eq-chp-1}
\limsup_{|z|\rightarrow1^{-}}\left\{(1-|z|^{2})\mbox{Re}\left(z\frac{h''(z)}{h'(z)}\right)\right\}\leq\sup_{z\in\mathbb{D}}\left\{(1-|z|^{2})\left|\frac{h''(z)}{h'(z)}\right|\right\}<2.
\ee It follows from (\ref{eq-chp-1}) and Theorem \ref{thm-4.0}(b)
that $\Omega$ is a radial John disk.

Now we prove the sharpness part. For $z\in\mathbb{D}$, let
$$f(z)=\frac{1}{2}\log\frac{1+z}{1-z}.$$ Then
$$\sup_{z\in\mathbb{D}}\left\{(1-|z|^{2})\left|\frac{f''(z)}{f'(z)}\right|\right\}=2,$$
and $f(\mathbb{D})$ is an infinite strip and hence not a radial John
disk. \hfill $\Box$




\subsection*{Conflict of Interests}
The authors declare that there is no conflict of interests regarding the publication of this paper.



\subsection*{Acknowledgements}
 This
research was partly supported by the Science and Technology Plan
Project of Hengyang City (No. 2018KJ125), the National Natural
Science Foundation of China (No. 11571216), the Science and
Technology Plan Project of Hunan Province (No. 2016TP1020), the
Science and Technology Plan Project of Hengyang City (No.
2017KJ183), and the Application-Oriented Characterized Disciplines,
Double First-Class University Project of Hunan Province
(Xiangjiaotong [2018]469).

\normalsize

\end{document}